\documentclass[leqno]{amsart}
\usepackage{amsmath,amsfonts,amssymb,amscd,amsthm,amsbsy,upref}
\newtheorem{thm}{Theorem}

\newtheorem{cor}[thm]{Corollary}
\newtheorem{prop}[thm]{Proposition}
\theoremstyle{definition}
\newtheorem{defin}[thm]{Definition}

\newtheorem{remark}[thm]{Remark}
\newtheorem{remarks}[thm]{Remarks}
\newtheorem{question}[thm]{Question}

\def\nat{{\mathbb N}}

\def\L{{\mathcal L}}

\catcode`@=11 \@mparswitchfalse  
\newcounter{mnotecount}[section]

\begin{document}
\title{On strongly asymptotic $\ell_p$ spaces and minimality}
\author{S. J. Dilworth}
\address{Department of Mathematics\\ University of South
Carolina,  Columbia, SC 29208.}
\email{$dilworth@math.sc.edu$}

\author{V. Ferenczi}
\address{Equipe d'Analyse Fonctionnelle, Bo\^\i te 186\\ Universit\'e Paris 6, 4, place Jussieu, 75252 Paris Cedex 05, France.}
\email{$ferenczi$@$ccr.jussieu.fr$}

\author{Denka Kutzarova}
\address{Institute of Mathematics, Bulgarian Academy of
Sciences, Sofia, Bulgaria.}
\curraddr{Department of Mathematics\\ University of Illinois at
Urbana-Champaign, Urbana, IL 61801.}
\email{$denka@math.uiuc.edu$}

\author{E. Odell}
\address{Department of Mathematics\\ The University of Texas at Austin\\
1 University Station C1200\\ Austin, TX 78712-0257.}
\email{$odell@math.utexas.edu$}

\thanks{Research of the fourth
named author was partially supported by the National Science Foundation.}

\begin{abstract}
Let $1 \le p \le \infty$ and let $X$ be a  Banach space with a strongly
asymptotic $\ell_p$ basis $(e_i)$. If
$X$ is minimal and $1 \leq p <2$, then $X$
is isomorphic to  a subspace of $\ell_p$. If $X$ is minimal  and $2
\leq p <\infty$, or if $X$ is complementably minimal and $1 \le p \le \infty$,
  then $(e_i)$ is equivalent to the unit vector basis of
$\ell_p$ (or $c_0$ if $p=\infty$).
\end{abstract}

\maketitle

\section{Introduction}

The notion of minimality was introduced by H.~Rosenthal. An
infinite-dimensional Banach
space $X$ is {\em minimal\/} if every infinite-dimensional subspace
has a further subspace isomorphic to $X$.

Let $1 \le p \le \infty$. 
A Banach space $X$ with a basis $(e_i)$ is {\em asymptotic\/} $\ell_p$
\cite{MT} if there exist $C <\infty$ and an increasing function 
$f:\nat \to \nat$ such that, for all $n \in \mathbb{N}$,  
every normalized block basis
$(x_i)_{i=1}^n$ of $(e_i)_{i=f(n)}^\infty$ is $C$-equivalent to the
unit vector basis of $\ell_p^n$. In this case $(e_i)$ is called an
{\em asymptotic\/} $\ell_p$ {\em basis\/} for $X$.

The only known examples of minimal spaces were $\ell_p$ ($1\le p
< \infty$) and $c_0$ and their subspaces until the original Tsirelson 
space $T^*$ \cite{CJT}, which is asymptotic $\ell_\infty$, was shown to 
be minimal \cite{CJT}. Tsirelson's space $T$ is
asymptotic $\ell_1$ and does not contain a minimal subspace \cite{CO}.

The next minimal space was constructed by Schlumprecht \cite{S}, and in
\cite{CKKM} a superreflexive version of $S$ was given. Both versions
are actually {\em complementably
minimal}, i.e. every infinite-dimensional subspace has a complemented
subspace isomorphic to the whole
space. Gowers
\cite{G} included minimality in his  partial classification of Banach
spaces, which motivated a series of results relating minimality to subsymmetry
\cite{P}, or to the number of nonisomorphic (resp. incomparable)
subspaces of a
Banach space \cite{FR},
  \cite{F} (resp. \cite{R1}).

We shall call a Banach space $X$ with a basis $(e_i)$
{\em strongly asymptotic\/} $\ell_p$
if there exist  $C <\infty$ and an increasing function
$f: \nat \to \nat$ such that, for all $n \in \mathbb{N}$, 
every normalized sequence
$(x_i)_{i=1}^n$ of disjointly supported vectors from
$[(e_i)_{i=f(n)}^\infty]$ is $C$-equivalent to the
unit vector basis of $\ell_p^n$.

In \cite{CO} it was proved that the standard bases of
$T$ and its strongly asymptotic version (sometimes called modified
Tsirelson's space) are equivalent. A new class of strongly
asymptotic $\ell_p$ spaces was introduced in \cite{ADKM}.
Sufficient conditions for selecting strongly asymptotic $\ell_p$ subspaces of
a given Banach space were
   given in \cite{FFKR} and \cite{T}.

Our main results are proved in \S3 (Theorem~\ref{pm}) and
\S4 (Theorem~\ref{cpm}). We summarize these results in slightly weaker form in

\begin{thm}\label{thm1}
Let $1 \le p \le \infty$ and let $X$ be a Banach space with a
strongly asymptotic $\ell_p$ basis $(e_i)$:
\begin{itemize}
\item if $X$ is minimal and  $1 \leq p <2$, then $X$ is
isomorphic to a subspace of $\ell_p$;
\item if $X$ is minimal  and $2 \leq p <\infty$, or if $X$ is
complementably minimal and $1\le p \le \infty$,
  then $(e_i)$ is equivalent to
the unit vector basis of
$\ell_p$ (or $c_0$ if $p=\infty$).
\end{itemize}
\end{thm}

For $1 \leq p <2$, we also give examples of  strongly asymptotic $\ell_p$
basic sequences in $\ell_p$ spanning minimal spaces that are not
isomorphic to $\ell_p$ (Theorem~\ref{thm: example}).

Recall again that  $T^*$ is a minimal space which is strongly asymptotic
$\ell_\infty$.

\begin{question}
Does there exist a  minimal asymptotic $\ell_p$ space that does 
not embed into $\ell_p$ for $1\le p <\infty$?
\end{question}

A broader definition of ``asymptotic $\ell_p$ space''
which is independent of a basis is given in \cite{MMT}.
The definition is that for some $C<\infty$ and
for all $n \in \mathbb{N}$, if $(e_i)_{1}^n \in \{X\}_n$, where
$\{X\}_n$ is the $n^{th}$
asymptotic class selected via the filter of  finite-codimensional
subspaces, then
$(e_i)_1^n$ is $C$-equivalent to the unit vector basis of $\ell_p^n$.
This notion does not seem to lend itself to a strongly asymptotic version.
But these spaces contain asymptotic $\ell_p$ basic sequences.

We use standard Banach space theory notation and terminology as in \cite{LT}.
All subspaces (of a Banach space $X$) are assumed to be closed and
infinite-dimensional unless stated otherwise. Let $(e_i)$ be a basis for $X$.
We say that a sequence $(x_k)$
of nonzero vectors is a \textit{block basis} if there exist
integers $0=n_0<n_1<\dots$ and scalars $(a_i)$
such that $$x_k= \sum_{i=n_{k-1}+1}^{n_k} a_i e_i\qquad (k\ge1).$$
The closed linear span of a block basis is called a \textit{block subspace}.

A minimal space is {\em C-minimal\/} if every subspace of $X$
contains a {\em C-isomorphic copy\/} of $X$, i.e. there exist $W \subseteq X$
and an isomorphism $T \colon X \rightarrow W$ with $\|T\|\|T^{-1}\| \le C$.
P.~Casazza proved that a minimal space must be $C$-minimal for some
$C<\infty$ \cite{C}.

The second author and C.~Rosendal \cite{FR} defined a space $X$ with a basis
$(e_i)$ to be {\em block minimal} (resp.\ {\em $C$-block minimal})
 if every block subspace of $X$ has a further block subspace isomorphic (resp.
$C$-isomorphic) to $X$;  also $X$ is defined to be {\em
equivalence block minimal} (resp {\em $C$-equivalence block
minimal}) if every block basis has a further block basis
equivalent (resp. $C$-equivalent) to $(e_i)$. They proved that
every block minimal (resp.\ equivalence block minimal) space is
$C$-block minimal (resp. $C$-equivalence block minimal) for some
$C < \infty$; they deduced that a basis which is asymptotic
$\ell_p$ for $1\le p \le \infty$ (in the broader sense associated
to the filter of tail subspaces) and equivalence block minimal
must be equivalent to the unit vector basis of $\ell_p$ or $c_0$.

We remark that $T^*$ is minimal but does not contain any block minimal
block subspace. Indeed, recall \cite{CS} that any block basis of $T^*$
is equivalent to a subsequence of the standard basis $(t_n^*)$ and
that the subspaces $[(t_{n_k}^*)]$ and $[(t_{m_k}^*)]$ are isomorphic
if and only if $(t_{n_k}^*)$ and $(t_{m_k}^*)$ are equivalent. 
Thus, if $T^*$ had a
block minimal block subspace $[(x_i^*)]$, then  $[(x_i^*)]$ would be
equivalence block minimal and therefore $C$-equivalence block minimal
for some $C<\infty$. That, combined with $(x_i^*)$ being an asymptotic
$\ell_\infty$ basis, would imply that $(x_i^*)$ is equivalent to the
unit vector basis of $c_0$, which is a contradiction since $T^*$ does not
contain $c_0$.

\section{Asymptotic $\ell_p$ subspaces of $L_p$}

Our reference for results concerning $L_p$ will mostly be the survey article by
D.~Alspach and the fourth author \cite{AO}. For the definition and properties
of the Haar basis $(h_i)$ we take \cite[Section~2.c]{LT2} as our reference.

We note that every strongly asymptotic
$\ell_p$ sequence $(e_i)$ is unconditional. Indeed, if $x$ and
$y$ are disjointly supported vectors belonging to $[(e_i)_{i \ge f(2)}]$, then
$$\|x \pm y\| \approx (\|x\|^p+\|y\|^p)^{1/p},$$ 
which implies unconditionality.

Let $X$ be a Banach space with a basis $(e_i)$. A \textit{blocking}
of $(e_i)$ is a finite-dimensional decomposition (FDD) for
$X$ corresponding to a partition $\{F_n\colon n \in \nat\}$ of $\nat$ into
successive finite subsets, i.e. $X=\sum_{n=1}^{\infty}[e_i\colon i \in F_n]$.
Let $T_{\omega}=\{(n_1,\ldots,n_k) \in {\nat}^k\colon k \in \nat\}$ be the
countably branching tree ordered by $(n_1,\ldots,n_k) \leq
(m_1,\ldots,m_l)$ if $k \leq l$ and $n_i = m_i$ for $i \leq k$.  We say
that $T \in T_{\omega}(X)$ if $T=\{x_{(n_1,\ldots,n_k)}:(n_1,\ldots,n_k)
\in T_{\omega}\} \subset S_X$.  $T$ is a
{\em block tree with respect to the basis\/} $(e_i)$
if in addition $(x_{(n)})_{n \in \nat}$ and 
$(x_{(n_1,\ldots,n_k,n)})_{n \in \nat}$
are all block bases of $(e_i)$ for all $(n_1,\ldots,n_k) \in T_{\omega}$.

\begin{prop}\label{aslp}
Let $1 < p <\infty$ and let $X$ be a subspace of $L_p[0,1]$ with an
asymptotic $\ell_p$ basis $(e_i)$. Then $X$
embeds into $\ell_p$. If $p \geq 2$ and $(e_i)$ is strongly
asymptotic $\ell_p$,
then $(e_i)$ is actually equivalent to the unit vector basis of $\ell_p$.
\end{prop}

\begin{proof} We may of course assume that $p \neq 2$. Then the space $X$
does not contain a subspace isomorphic to $\ell_2$, otherwise some block
basis of $(e_i)$ would be equivalent to the unit vector basis of $\ell_2$,
contradicting the fact that $(e_i)$ is asymptotic $\ell_p$.
If $p>2$ this implies by \cite[Theorem~30]{AO} that
$X$ embeds into $\ell_p$ (in fact, $X$ $1+\varepsilon$-embeds into $\ell_p$ for all $\varepsilon>0$
\cite{KW}).
Assume now that $1 < p <2$. 
By a theorem of W.~B.~Johnson \cite[Theorem~30]{AO}
it suffices to prove that for some $C<\infty$ every normalized 
block basis $(x_i)$ of $(e_i)$ admits a subsequence $C$-equivalent 
to the unit vector basis of $\ell_p$.
Passing to a subsequence we may assume that $(x_i)$ is a 
perturbation of a block basis of the Haar basis and
hence is $C(p)$-unconditional for some constant $C(p)$ depending 
only on $p$ \cite[Theorem~2.c.5]{LT2}.
By \cite[Lemma~28]{AO} it suffices to show that there exist $\delta(K)>0$,
where $K$ is the asymptotic $\ell_p$-constant of $(e_i)$, 
$(y_i) \subseteq (x_i)$ and disjoint measurable
sets $E_i$ in $[0,1]$ with $\|y_{i | E_i}\| \ge \delta(K)$ for all $i$.

By a theorem of Dor \cite{D} there exists $\delta(K)>0$ such that if $(x_i)_{i \in A}$ is $K$-equivalent to the unit vector basis of
$\ell_p^{|A|}$ then there exist disjoint measurable sets $(F_i)_{i \in A}$ in $[0,1]$ with $\|x_{i | F_i}\| > \delta(K)$ for
$i \in A$. Since $(e_i)$ is $K$-asymptotic $\ell_p$ we can thus, for all $n$ and $\varepsilon>0$,  choose $i >n$ and a
measurable set $F_i$ with $\lambda(F_i)<\varepsilon$ (where $\lambda$ denotes Lebesgue measure) and $\|x_{i|F_i}\|>\delta(K)$. We then proceed inductively.
Let $y_1=x_1$ and $F^1_1 = [0,1]$. Assume that $(y_1,\dots,y_n)$ have been  chosen along with disjoint sets
$F^n_1,\dots, F^n_n$ so that $\|y_{i|F^{n}_i}\| > \delta(K)$. Applying the above for a suitable $\varepsilon$, using the uniform integrability
of $(|y_i|^p)_{i=1}^n$, we can choose $y_{n+1}$ and $F^{n+1}_{n+1}$ so that
$\|y_{n+1|F^{n+1}_{n+1}}\|>\delta(K)$, $(y_i)_1^{n+1}$ is a subsequence of $(x_i)$, and $\|y_{i|F^{n+1}_i}\| > \delta(K)$ for $i \le n$,
where $F^{n+1}_i = F^{n}_i \setminus F^{n+1}_{n+1}$. The claim now follows by setting $E_i = \cap_{n \ge i} F^n_i$.

Finally, assume that $p>2$ and that $(e_i)$ is strongly asymptotic $\ell_p$ with constant
$C<\infty$ and an associated function $f$. As we noted, $(e_i)$ is
unconditional. By
\cite{DS} (see also \cite{JMST} or \cite[Corollary~2.e.19]{LT2} for a more
general lattice version) a normalized unconditional sequence $(x_i)$ in
$L_p$ is either equivalent to the unit vector basis of $l_p$, or it
has, for all $n \in \nat$, $n$ disjoint normalized blocks which form
a sequence $2$-equivalent to the unit vector basis of $l_2^n$.

Let $k \in \nat$ be such that the unit vector bases of $l_2^k$ and
$l_p^k$ are not
$2C$-equivalent. Since $(e_i)$ is strongly
asymptotic $\ell_p$, there cannot exist a normalized length $k$
sequence of disjoint
blocks on $(e_i)_{i \geq f(k)}$ which is $2$-equivalent to the basis
of $l_2^k$ . So
$(e_i)_{i \geq f(k)}$, and therefore $(e_i)$, is equivalent to the
unit vector basis
of
$l_p$.
\end{proof}

Note that it remains open whether we may remove ``strongly'' from the
last part of Proposition \ref{aslp}. \begin{question}
 Does there exist an unconditional  asymptotic $\ell_p$ basic sequence in
$L_p$, $p  > 2$, which is not equivalent to the
unit vector basis of $\ell_p$? \end{question}

\begin{prop} \label{Lonecase} 
Let $X$ be a subspace of $L_1[0,1]$. If $X$ has
an unconditional asymptotic $\ell_1$ basis then $X$ embeds into $\ell_1$.
\end{prop}

\begin{proof} Let $(e_i)$ be an unconditional asymptotic $\ell_1$ basis for $X$.
Then $(e_i)$ is necessarily boundedly complete, so by \cite[Proposition~31]{AO}
it suffices to prove that for some $C<\infty$ if $T \in T_\omega(X)$ is
a block tree with respect to $(e_i)$ then some branch of $T$ is $C$-equivalent
to the unit vector basis of $\ell_1$. The proof of this is much like the
previous proof in the case $1<p<2$. We need only show that there exists
 $\delta=\delta(K)>0$ so that for all such $T$ we can find a branch $(y_i)$
and disjoint measurable sets $(E_i)$ with $\|y_{i|E_i}\| > \delta$
for all $i$. Again we use Dor's theorem. If $(y_i)_1^n$ have been
chosen as an initial segment in $T$ along with disjoint sets
$(F^n_i)_{i=1}^n$ satisfying  $\|y_{i|F^n_i}\| > \delta$ for $i
\le n$ then we can select $y_{n+1}$ from among the successors of
$y_n$ and $F^{n+1}_{n+1}$ with $\|y_{n+1|F^{n+1}_{n+1}}\| >
\delta$ and $\lambda(F^{n+1}_{n+1})$ as small as we please. In
particular we may insure that $\|y_{i|F^n_i \setminus
F^{n+1}_{n+1}}\| > \delta$ for $i \le n$. Then as earlier we let
$E_i = \cap_{n \ge i} F^n_i$.
\end{proof}

\begin{remarks}  (1) We do not know if Proposition~\ref{Lonecase} holds if we merely assume that $(e_i)$ is asymptotic $\ell_1$.
Then $(e_i)$ would still be boundedly complete but we need unconditionality to conclude that the branch $(y_i)$
with $\|y_{i|E_i}\|>\delta$ is $C(K)$-equivalent to the unit vector basis of $\ell_1$ using \cite[Lemma~28]{AO}.

However, we can drop ``unconditionality'' in any of the following three cases:
\begin{itemize}
\item[(a)]  $(e_i)$ is a block basis of $(h_i)$ or more generally for all $n \in \mathbb{N}$
$$ \lim_{m \rightarrow \infty} \sup \{ \|P_n x\| \colon x \in B_{[ e_i]_{i=m}^\infty} \}=0,$$
where $P_n$ is
the basis projection of $L_1[0,1]$ onto the span of the first $n$ Haar functions;
\item[(b)] for some $K<2$ every normalized block basis of $(e_i)$ admits a subsequence
$K$-equivalent to the unit vector basis of $\ell_1$.
\item[(c)] $(e_i)$ is $K$-asymptotic $\ell_1$ with $K < 2$.
\end{itemize}

Indeed, the branch $(y_i)$ obtained in the proof can be chosen so that \begin{align*}
\|\sum_i a_i y_i \|\ge \| \sum_i a_i y_{i|_{\cup_jE_j}}\|
\gtrsim \sum_i \|a_i y_{i|E_i} + \sum_{j>i} a_j y_{j |E_i}\|. \end{align*}
Now under (a) we can also choose $(y_i)$ so that, for all $i$,  $(y_{j|E_i})_{j \ge i}$ is $2$-basic
and so
$$\|a_i y_{i|E_i} + \sum_{j>i} a_j y_{j |E_i}\| \ge \frac{1}{2} |a_i| \|y_{i|E_i}\| \ge \frac{1}{2} |a_i| \delta.$$
Next we shall show that case (b) yields that we can choose the branch $(y_i)$ and sets $(E_i)$ so that $\|y_{i|E_i}\|>\delta>1/2$ if
$1/2 < \delta <1/K$. Since \begin{align*}
\sum_i\|a_i y_{i|E_i} + \sum_{j>i} a_j y_{j |E_i}\| &\ge \sum_i (|a_i| \delta - \sum_{j>i} |a_j| \|y_{j|E_i}\|)\\
& \ge \delta \sum_i  |a_i| - \sum_{i} |a_i| \|y_{i|[0,1] \setminus E_i}\| \\
& \ge (\delta - (1-\delta)) \sum_i |a_i| = (2\delta-1) \sum_i |a_i|, \end{align*}
we obtain that $(y_i)$ is equivalent to the unit vector basis of $\ell_1$.

Suppose that $(x_i)$ is a normalized block basis of $(e_i)$ and (b) holds. Passing to a subsequence we may assume that $(x_i)$ is $K$-equivalent to
the unit vector basis of $\ell_1$. Moreover by
Rosenthal's ``subsequence splitting lemma'' (see e.g.\  \cite{GD} or \cite[p. 135]{JMST} for a general lattice version)
 we may assume that $(x_i) = (u_i + z_i)$
where $(z_i) = (x_{i|E_i})$ is disjointly supported, $\lim_i \|z_i\| = b$, and $(|u_i|)$ is uniformly integrable. It suffices
to prove that $b \ge 1/K$. For then the proof of Proposition~\ref{Lonecase} will yield the $\|y_{i|E_i}\|>\delta$ claim.
Given $\varepsilon>0$ we can form an absolute convex combination $\sum_i a_i x_i$ with $\|\sum_i a_i u_i\| < \varepsilon$, and thus
$\|\sum_i a_i z_i\| = \sum_i |a_i|\|z_i\| \ge 1/K - \varepsilon$. Since these can be found arbitrarily far out we deduce that $b \ge 1/K$.

Finally assume (c) holds.
Let $(x_i)$ again  be a normalized block basis of $(e_i)$. We may assume, passing to a subsequence, as in case (b)
that $(x_i) = (u_i + z_i)$,
$(z_i) = (x_{i|E_i})$ is disjointly supported, $\lim_i \|z_i\| = b$,  $(|u_i|)$ is uniformly integrable, and moreover $u_i = u+v_i$ where
$(v_i)$ is a weakly null perturbation of a block basis of the Haar basis for $L_1[0,1]$. We assume  $\|v_i\| \rightarrow a \in [0,2]$
and first assume  $a \ne 0$. Passing again to a subsequence we may assume $(v_i)$ has a spreading model $(\tilde v_i)$
which is unconditional, subsymmetric, and (e.g.\ by Dor's theorem \cite{D}) is not equivalent to the unit vector basis of $\ell_1$.
Thus $\lim_n (1/n)\|\sum_{i=1}^n \tilde v_i\| = 0$ (see e.g.\ \cite{BL}). In particular, we can thus find (even if $a=0$) for all
$\varepsilon>0$ an absolute convex combination $\sum_{i \in E_\varepsilon} a^\varepsilon_i x_i$ with $\|\sum_{E_\varepsilon} a^\varepsilon_i u_i\| < \varepsilon$ and
$f(|E_\varepsilon|) \le \min E_\varepsilon$,
 where $f$ is the function associated to the asymptotic $\ell_1$ sequence $(e_i)$. As in case (b) we deduce that $b \ge 1/K$ and finish the proof as in case (b).

(2) Note that if $(e_i)$ is a $K$-asymptotic $\ell_1$ basis for  $X \subset L_1[0,1]$ then $X$ has the
strong Schur property. To see this, let $\delta>0$ and let $(x_i)$ be a sequence in $B_X$ satisfying $\|x_i - x_j\| \ge \delta$.
Since $(e_i)$ is boundedly complete, we may assume (by passing to a subsequence) that   $x_i = y + z_i$, where $y, z_i \in X$
and $(z_i)$ is a perturbation of a block basis of $(e_i)$.
We may assume that $\|z_i\| \rightarrow a \ge \delta/2$ and that $(z_i/\|z_i\|)$ is $K'$-asymptotic $\ell_1$ where $K'>K$ is arbitrary.
As in (1) we may also assume, by passing to a subsequence, that there exist disjoint sets $(E_i)$ with
$\lim_i \|(z_i/\|z_i\|)_{|E_i}\| = b \ge 1/K'$ and so $\lim_i \|z_{i|E_i}\| = ab \ge \delta/(2K')$.
By Rosenthal's lemma \cite{R} some subsequence of $(z_i)$ is $(4K'/\delta)$-equivalent to the unit vector basis of $\ell_1$.
Hence a subsequence of $(x_i)$ is $(4K'(B+2)/\delta)$-equivalent to the unit vector basis of $\ell_1$, where $B$ is the basis constant of
$(e_i)$.

In \cite{BR} some $1$-strong Schur subspaces of $L_1[0,1]$ which do not embed into $\ell_1$
are constructed.

(3) Let $(E_i)$ be an FDD for a Banach space $X$.
We say that  $(E_i)$ is an asymptotic $\ell_p$ (resp.\ strongly asymptotic $\ell_p$) FDD
 if there exist $C <\infty$ and an increasing function
$f: \nat \to \nat$ such that, for all $n \in \mathbb{N}$, every normalized block basis (resp.\ sequence of disjointly supported vectors)
$(x_i)_{i=1}^n$ of $[(E_i)_{i=f(n)}^\infty]$ is $C$-equivalent to the
unit vector basis of $\ell_p^n$. The proof of Proposition~\ref{aslp} carries over to show that if $X$ is a subspace of $L_p$  ($1<p<\infty$) with an
asymptotic $\ell_p$ FDD then $X$ embeds into $\ell_p$. Similarly, Proposition~\ref{Lonecase} remains valid if ``basis'' is replaced by ``FDD''.
\end{remarks}

\section{Strongly asymptotic $\ell_p$ spaces and primary minimality}

A key ingredient of the proof of Theorem~\ref{pm} is a uniformity lemma about
embedding into subspaces generated by tail subsequences of a basis
$(e_i)_{i \in \nat}$
(i.e. subsequences of the form $(e_i)_{i \geq k}$, for $k \in \nat$).

We recall that a Banach space $X$ is said to be {\em primary} if whenever
$X=Y \oplus Z$,
then either $X$ is isomorphic to $Y$ or $X$ is isomorphic to $Z$.

\begin{defin} \label{primarilyminimal}
A Schauder basis $(e_n)_{n \in \nat}$ is said to be {\em primarily
minimal} if whenever
$(I,J)$ is a partition of $\nat$, then either $[e_n]_{n \in \nat}$ embeds into
$[e_n]_{n \in I}$
or $[e_n]_{n \in \nat}$ embeds  into $[e_n]_{n \in J}$. A Banach space $X$ is {\em primarily
minimal} if whenever $X=Y
\oplus Z$, then either $X$ embeds into $Y$ or $X$ embeds into $Z$.
\end{defin}

Obviously any unconditional basis of a primarily minimal space must
be primarily minimal.
Also the class of primarily minimal spaces contains the class of
primary spaces as well as
the class of minimal spaces.

\begin{prop}\label{prop}
Let $X$ be a Banach space with a primarily minimal
basis $(e_n)_{n \in \nat}$. Then there exists
$C<\infty$ such that $X$ $C$-embeds into all of its tail subspaces.
\end{prop}

\begin{proof}
Let $2^{\omega}$ be equipped with its usual topology, i.e. basic open sets
are of the form $N(u)=\{\alpha=(\alpha_k)_{k \in \nat} \in
2^{\omega}: \forall k \leq m, \alpha_k=u_k \}$, where
$u=(u_k)_{k
\leq m}$ is a  sequence of $0$'s and $1$'s of length $m$.

Let $A$ be the set of $\alpha$'s in $2^{\omega}$ such that $X$ embeds into
$[e_i]_{\alpha_i=1}$, and let $A'$ be the set of $\alpha$'s in $2^{\omega}$ such
that $X$ embeds into $[e_i]_{\alpha_i=0}$.
We also let, for $n \in \nat$, $A_n$ (resp. $A^{\prime}_n$) be the
set of $\alpha$'s such that $X$ $n$-embeds into
$[e_i]_{\alpha_i=1}$ (resp.\ $[e_i]_{\alpha_i=0}$).

Now $(e_n)$ being primarily minimal means that $2^{\omega}=A \cup A'=
(\cup_{n \in \nat}A_n)\cup(\cup_{n \in \nat}A^{\prime}_n)$. By the
Baire Category Theorem
we deduce that the closure of  $A_N$ or $A^{\prime}_N$ has nonempty
interior for some $N
\in \nat$, and hence $A_N$ or $A^{\prime}_N$ is dense in some basic
neighborhood.
We note that the map $f:2^{\omega} \rightarrow
2^{\omega}$ defined by $f((\alpha_i)_{i \in \nat})=(1-\alpha_i)_{i
\in \nat}$ is a
homeomorphism such that
$f(A_N)=A^{\prime}_N$. It follows that  $A_N$
is dense in some
basic neighborhood. Let $u=(u_i)_{i \leq m}$ be such that $A_N$ is dense
in $N(u)$.

Let $E = \{i : u_i = 1, i \le m \}$. It follows that $X$ $N$-embeds
into $[(e_i)_{i \in E} \cup (e_i)_{i > k}]$ for all $k \in
\nat$. Using the uniform equivalence between $(e_i)_{i \in E}$ and
$(e_i)_{i \in F}$ for all $|F| = |E|$, we see that there exists $C$
such that $X$ $C$-embeds into $[(e_i)_{i > k}]$ for all $k$.
\end{proof}

\begin{thm}\label{pm} 
Let $X$ have a primarily minimal, strongly asymptotic $\ell_p$ 
basis $(e_i)$ for some $1 \leq p <\infty$. If $1 \leq p<2$ then $X$
embeds into $\ell_p$,  and if $p \geq 2$ then $(e_i)$ is
equivalent to the unit vector basis of $\ell_p$.
\end{thm}

\begin{proof}

We may assume that $X$ has a 1-unconditional basis $(e_i)$
with a strongly asymptotic $\ell_p$ constant $C$.

By Proposition \ref{prop} there exists $K<\infty$  such that $X$
$K$-embeds into each of  its tail subspaces
associated to $(e_n)_{n \in \nat}$.
By a result of Johnson \cite{J}, for each $n$
  there exists  $N(n)$ such that for any
$n$-dimensional subspace $F$ of any Banach space
$E$ with a $1$-unconditional basis there exist $N(n)$ normalized
disjointly supported vectors $(w_i)_{i=1}^{N(n)}$ (with respect to
the basis) such that
$F$ is 2-isomorphic to a subspace of $[(w_i)_{i=1}^{N(n)}]$. Let
$F$ be an arbitrary finite-dimen\-sional subspace of $X$ and let $n$ be
its dimension. Consider $Y=[(e_i)_{i \ge f(N(n))}]$. Then $Y$
contains a $K$-isomorphic copy of $X$, and thus also of $F$. By the
above result, there exist $N=N(n)$ normalized disjointly supported
vectors $(w_i)_{i=1}^N$ of $[(e_i)_{i \ge f(N)}]$ such that $F$ is
$2K$-isomorphic to a subspace of $W:=[(w_i)_{i=1}^N]$. Therefore, $W$ is
$C$-isomorphic to $\ell_p^N$. Thus, $F$ is $2CK$-isomorphic to
a subspace of $\ell_p^N$. It follows that $X$ is crudely finitely representable
in $\ell_p$ and hence embeds isomorphically into $L_p[0,1]$ \cite{LP}.
The result now follows from Proposition~\ref{aslp} and 
Proposition~\ref{Lonecase}.
\end{proof}

\begin{remark}\label{remark} For $1 \leq p <2$, we actually obtain that
$(e_i)$ may be blocked into a decomposition
$\sum_{i=1}^{\infty}\oplus F_i \cong  (\sum_{i=1}^{\infty}\oplus F_i)_p$,
where the $F_i$'s embed uniformly into
$\ell_p$.
Indeed,  any unconditional basis for a subspace $X$ of
$\ell_p$, $1 \leq p<\infty$, may be blocked into such an FDD. For
$p>1$ this is due to
W. B. Johnson and M. Zippin \cite{JZ}, and unconditionality is not
required. For $p=1$, this follows for example from
\cite[Proposition~31]{AO}.
\end{remark}

As we mentioned in the Introduction, we do not know if a
positive result similar to Theorem \ref{thm1} is true for the more
general class of asymptotic
$\ell_p$ spaces, $1\le p < \infty$. For example, we do not know what
the case is with the Argyros-Deliyanni mixed Tsirelson space $X_u$
\cite{AD}. Two main ingredients were used in \cite{CJT} to prove the
minimality of $T^*$, namely the universality of $\ell_\infty^n$'s for
all finite dimensional spaces and an appropriate blocking
principle. For $X_u$ there is no corresponding blocking principle.
On the other hand, it was proved in \cite{ADKM} that all of its 
subspaces contain
uniformly $\ell_\infty^n$'s, which makes it impossible to use a
local argument as in the case of strongly asymptotic $\ell_p$ spaces.

Note that there is no version  of Theorem~\ref{pm} for spaces with a 
primarily minimal FDD (which is
defined by replacing the basis $(e_i)$ by an FDD $(E_i)$ throughout 
Definition~\ref{primarilyminimal}). 
Indeed, let $X = (\sum_{n=1}^\infty \oplus \ell_\infty^n)_p$, 
where $1 \le p < \infty$. 
Then the natural FDD for  $X$ is easily seen to be a primarily
minimal strongly asymptotic $\ell_p$  FDD, but obviously  $X$ does not 
embed into $\ell_p$. However, we have the following result for
minimal spaces.

\begin{prop}
Suppose that $X$ is minimal and that $(E_i)$ is  a strongly 
asymptotic $\ell_p$ FDD for $X$, where $1 \le p < \infty$.
Then $X$ embeds into $\ell_p$. 
\end{prop}

\begin{proof} 
It suffices to show that $\ell_p$ embeds into $X$.
Choose $e_i \in E_i$ for $i \ge 1$ with $\|e_i||=1$. Then $(e_i)$ is a strongly asymptotic $\ell_p$ sequence spannning a
minimal subspace $Y$ of $X$, so $Y$ embeds into $\ell_p$  by Theorem~\ref{pm}, and hence $\ell_p$ embeds into $X$.
\end{proof}

We now present examples of minimal strongly asymptotic $\ell_p$
spaces that are not
isomorphic to $\ell_p$, for $1 \leq p<2$. Note that if $(e_i)$ is a strongly
asymptotic
$\ell_p$ basis ($1\le p \le \infty$), then $(e_i^*)$ is a strongly
asymptotic $\ell_q$
($1/p+1/q=1$)  basic sequence with the same $f:\nat \rightarrow \nat$.
  This follows easily from the unconditionality of $(e_i)$ and a
standard H\"older's
inequality calculation.

\begin{thm} \label{thm: example} 
Let $1 \le p < 2$. There exists a minimal Banach space $X$ with
a strongly asymptotic $\ell_p$ basis $(e_i)$ satisfying the following: 
\begin{itemize}
\item[(i)] $X$ is not isomorphic to $\ell_p$;\\
\item[(ii)] If $p>1$ then $X^*$ does not embed into $L_q$ ($1/p+1/q=1$). 
\end{itemize}
\end{thm}

\begin{proof} 
We shall use the following two facts (the first
  follows from \cite{KP}, and the second follows from the existence, due
to Paul L\'evy, of $s$-stable random variables for $1 < s <2$): firstly,
if $s \notin \{2,p\}$ and $M>0$ then there exists $N$ such that $\ell_p$
does not contain an $M$-complemented
$M$-isomorphic copy of $\ell_s^N$;
secondly, if $p<s<2$ then $\ell_p$ contains almost isometric copies
(i.e. $(1+\varepsilon)$-isomorphic copies for all $\varepsilon>0$) of
$\ell_s^n$ for all $n$.  For each $n$,
let $s_n$ be defined by the equation
$1/s_n:=1/p-(1/p-1/2)/(2n)$. Note that $p<s_n<2$ and that $s_n
\rightarrow p$ rapidly
enough to ensure that the standard bases of
$\ell_{s_n}^n$ and $\ell_p^n$ are $2$-equivalent (as is easily
checked). By the first
fact applied to $M=n$ and $s=s_n$,  there exists $N(n)$ such that $X_n
:=\ell_{s_n}^{N(n)}$ is not
$n$-isomorphic to an $n$-complemented subspace of $\ell_p$.
Let $X := (\sum_{n=1}^\infty \oplus X_n)_p$ and let $(e_i)$ be the
basis of $X$ obtained by
concatenation of the standard bases of  $X_1, X_2,\dots$ in that order.

  By the second fact each $X_n$ is almost isometric to a subspace of
$\ell_p$, so $X$ is also almost isometric
to a subspace of $\ell_p$. Since $X_n$ is $1$-complemented in $X$ but
is not $n$-isomorphic to an
$n$-complemented subspace of $\ell_p$, it follows that $X$ is not
isomorphic to $\ell_p$, which proves (i).
  Property (ii) will follow from (i) and the fact that $(e_i)$ is a
strongly asymptotic $\ell_p$ basis: then
$X^*$ is a strongly asymptotic $\ell_q$ space which is not isomorphic
to $\ell_q$, and therefore
does not embed into $L_q$ by Proposition \ref{aslp}.

Finally, we verify that $(e_i)$ is a strongly asymptotic $\ell_p$
basis. Suppose that $n \ge
1$  and that the vectors
$x_1,\dots, x_n$ are disjointly supported vectors (with respect to
$(e_i)$) which belong to the tail space $(\sum_{j=n}^\infty \oplus
X_j)_p$. Write $x_i = \sum_{j=n}^\infty x^j_i$, where $x^j_i \in X_j$.
Then
\begin{align*}
\|\sum_{i=1}^n x_i\|^p &= \sum_{j=n}^\infty \|\sum_{i=1}^n x^j_i\|^p\\
&= \sum_{j=n}^\infty (\sum_{i=1}^n \|x^j_i\|^{s_j})^{p/s_j}\\
\intertext{(by the disjointness of $x^j_1,\dots, x^j_n$)}
&\approx \sum_{j=n}^\infty \sum_{i=1}^n \|x^j_i\|^p\\
\intertext{(by the $2$-equivalence of the standard bases
of $\ell_p^n$ and $\ell_{s_j}^n$
for all $j \ge n$)}
&= \sum_{i=1}^n \sum_{j=n}^\infty \|x^j_i\|^p\\
&= \sum_{i=1}^n \|x_i\|^p.
\end{align*}.
\end{proof}

\section{Strongly asymptotic $\ell_p$-spaces and complementation}

Let us define a basis $(e_n)_{n \in \nat}$ to be {\em primarily
complementably minimal} if
for any partition $(I,J)$ of $\nat$, $[e_n]_{n \in \nat}$ is
isomorphic to a complemented subspace
of $[e_n]_{n \in I}$ or
of $[e_n]_{n \in J}$. A space $X$ is {\em primarily complementably
minimal} if
whenever $X=Y \oplus Z$, then either $X$ embeds complementably into $Y$ or
$X$ embeds complementably into  $Z$.

Every Banach space which is primary, or which is complementably
minimal, is primarily complementably
minimal. We now prove a theorem which yields the complementably
minimal part of Theorem \ref{thm1}.

Recall that an unconditional basis $(x_n)_{n \in \nat}$ of a Banach space is
{\em sufficiently lattice-euclidean} if it has, for some $C \geq 1$
and every $n \in \nat$,
a $C$-complemented, $C$-isomorphic copy of $\ell_2^n$
whose basis is disjointly supported on $(x_n)_{n \in \nat}$. See
\cite{CK} for a general definition in the
lattice setting.
Note that strongly asymptotic $\ell_p$ bases for $p \neq 2$ are not
sufficiently lattice-euclidean.

\begin{prop}\label{propcomp}
Let $X$ be a Banach space with an unconditional, primarily
complementably minimal
basis $(e_n)_{n \in \nat}$. Then there exists
$K<\infty$ such that $X$ $K$-embeds as a $K$-complemented subspace
of its tail subspaces.
\end{prop}

We skip the proof since it is exactly the same as the proof of
Proposition~\ref{prop}, mutatis mutandis.
Note that unconditionality is required to preserve complemented
embeddings with uniform constants in the end
of the argument.

\begin{thm}\label{cpm} Let $X$ be a primarily complementably minimal
Banach space with a strongly asymptotic
$\ell_p$ basis $(e_n)$, $1 \leq p \leq \infty$. Then $(e_n)$ is
equivalent to the unit vector basis of
$\ell_p$ (or $c_0$ if $p=\infty$).
\end{thm}

\begin{proof} We may assume that $p \neq 2$ and we may
 also assume  that $(e_n)$ is
$1$-unconditional. By  Proposition \ref{propcomp}, there exists
 $K<\infty$ such that $X$ is $K$-isomorphic to some
$K$-complemented subspace of any tail
subspace $Y$ of $X$. Since $p \neq 2$, $X$ is
not sufficiently lattice euclidean; note
  also that the canonical basis of $Y$ is $1$-unconditional. Therefore
we deduce from
\cite[Theorem~3.6]{CK}  that $(e_n)$ is $c$-equivalent to a sequence
of disjointly supported vectors in
$Y^N$, for some
$c$ and $N$ depending only on $X$ and $K$.
Here
 $Y^N$ is equipped
with the norm $\|(y_i)_{i=1}^N\| = \max_{1\le i \le N} \|y_i\|$
and with the canonical basis obtained from the basis $(b_i)$ of
$Y$ with the ordering $(b_1,0,\dots,0), (0,b_1,0\dots,0),\dots,
(0,\dots,0,b_1), (b_2,0,\dots,0)$, etc. It is easy to check that
$Y^N$ is also strongly asymptotic $\ell_p$ with constant $C$ and
function $f$ depending only  on $N$ and on the strongly asymptotic
$\ell_p$ constant and function of $(e_n)$.

  Let $k \in \nat$ be arbitrary and let $Y$ be a tail subspace of $X$
such that $Y^N$ is
supported after the $f(k)$-th vector of the basis of $X^N$. The
sequence $(e_i)_{i \leq k}$ is therefore
$c$-equivalent to a
  sequence which is disjointly supported after $f(k)$, so is
$cC$-equivalent to the unit vector basis
of $l_p^k$.
\end{proof}

\begin{remark}  We applied  \cite[Theorem~3.6]{CK} with  $E_n = \nat$ for all $n$. There is an
 inaccuracy
 in
the statement of the theorem: the sets $E_n$ do not need to be assumed
disjoint. Theorem~3.6 is based on  \cite[Lemma~3.3]{CK}, and we want to point
out an imprecision in the proof of Lemma~3.3. The inequality on p.150,
line 11, goes the opposite way. In the special case in which we apply Theorem 3.6, however, we have
 $h_n=|Se_n||T^*e_n|$ and the proof is clear.
Actually, as Nigel Kalton has informed us, Lemma 3.3 is true as stated, but one needs
a small adjustment in the general case. More precisely, one should define
the functions $f_n$ and $g_n$ in a slightly different way. If
$0 \le h_n \le |Se_n||T^*e_n|$, we have that
$h_n$ can be expessed in the form
$h_n =f_n g_n$, where $|f_n|\le |Se_n|$ and $|g_n|\le |T^*e_n|$.
Then the proof is correct as it stands. \end{remark}

\section{Some consequences about the number of nonisomorphic
subspaces of a Banach space}

The question of the number of mutually nonisomorphic subspaces of a
Banach space
which is not isomorphic to
$\ell_2$ was first investigated by the second author and Rosendal
\cite{FR1},
\cite{FR}, \cite{R1} (but some ideas originated from an earlier paper
of N. J. Kalton \cite{K}). A separable
Banach space
$X$ is said to be {\em ergodic} if the relation $E_0$ of eventual
agreement between
sequences of
$0$'s and
$1$'s is Borel reducible to isomorphism between subspaces of $X$;
this means that there
exists a Borel map $f$ mapping elements of $2^{\omega}$ to subspaces
of $X$ such that
$\alpha E_0 \beta$ if and only if $f(\alpha) \simeq f(\beta)$. We refer to
\cite{FR} or \cite{R1} for
detailed definitions. We just note here that
an ergodic Banach space $X$ must contain $2^{\omega}$ mutually
nonisomorphic subspaces, and
furthermore, that it admits no Borel classification of isomorphism
classes by real numbers,
i.e. no Borel map
$f$  mapping subspaces of
$X$ to reals, with
$Y \simeq Z$ if and only if $f(Y)=f(Z)$. A natural conjecture is to ask if any
Banach space
nonisomorphic to $\ell_2$ must be ergodic.

\begin{thm} Let $1 \leq p \leq \infty$ and let $X$ be a Banach space
with a strongly
asymptotic
$\ell_p$ basis $(e_i)$. Then $(e_i)$ is equivalent to
the unit vector basis of $\ell_p$ (or $c_0$ if $p=\infty$), or
$E_0$ is Borel reducible to  isomorphism between subspaces of $X$ spanned
by subsequences of the basis (and in particular there are continuum many
mutually nonisomorphic complemented subspaces of $X$).
\end{thm}

\begin{proof}
The basis $(e_n)$ is unconditional. If $E_0$ is not Borel reducible to
isomorphism between subspaces of $X$ spanned by subsequences of the
basis, then by
\cite{FR1}, \cite{R1}, there exists $K<\infty$ such that the set
$\{\alpha: [e_n]_{\alpha_n=1} \simeq^K X\}$ is comeager in $2^{\omega}$. In
particular the set $A_K$ of $\alpha$'s such that $X$ $K$-embeds
$K$-complementably in $[e_n]_{\alpha_n=1}$
is comeager, thus dense, and the proof of Proposition \ref{propcomp}  applies
to deduce that for some $K'<\infty$, $X$ is $K'$-isomorphic to some
$K'$-complemented subspace of any tail
subspace of $X$. Then the proof of Theorem \ref{cpm}, if $p \neq 2$,
or Theorem \ref{pm}, if
  $2 \leq p < \infty$, applies.
\end{proof}

A consequence of this result is that the versions $T_p$ of
Tsirelson's space are
ergodic for $1 < p<\infty$. For $T$, a stronger result was already
proved by Rosendal
\cite{R0}.
Another consequence is that the mixed Tsirelson spaces and their
$p$-convexifications for $1 < p<\infty$ \cite{ADKM} are also ergodic.
For a space $X$ with a strongly asymptotic $\ell_p$ FDD $(F_i)$, we deduce from
Theorem 16 that $X$ is ergodic, or that $X=\sum \oplus F_i \simeq (\sum \oplus
F_i)_p$
in which case $X
\simeq l_p(X)$ by \cite{FG2} Corollary 2.12 (with the usual modifications when
$p=\infty$).

\begin{cor} Let $1 \leq p \leq\infty$ and let $X$ be a Banach space
with a strongly asymptotic
$\ell_p$ basis. Then
$X$ is isomorphic to $\ell_2$ or $X$ contains $\omega_1$
nonisomorphic subspaces.
\end{cor}

\begin{proof} Assume $X$ contains no more than countably many
mutually nonisomorphic
subspaces. By the above, $X$  is
  isomorphic to
$\ell_p$ (or $c_0$ if $p=\infty$). It is known that for $p \neq 2$,
$\ell_p$ contains at least $\omega_1$
nonisomorphic subspaces
\cite{LT, FG1} (in fact, $c_0$ and  $\ell_p$, $1
\leq p<2$, are ergodic). So
  $X$ is isomorphic to $\ell_2$.
\end{proof}

\begin{cor} Let $1 \leq p \leq \infty$. Let $X$ be a Banach space with a
strongly asymptotic
$\ell_p$
FDD. Then
$X$ is isomorphic to $\ell_2$ or $X$ contains $\omega$ nonisomorphic subspaces.
\end{cor}

\begin{proof} By the above we may assume that $X \simeq \ell_2(X)$. If $X$ has
finite cotype and $X$ is not
isomorphic to $\ell_2$,
then by \cite{A}, $\ell_2(X)$ and therefore $X$ contains at least $\omega$
nonisomorphic subspaces.
If $X$ does not have finite cotype, it contains $\ell_{\infty}^n$'s uniformly
and therefore
$X$ contains copies of the space $Y_p=(\sum_{n \in \nat}\ell_p^n)_2$ for any $p
\in [1,\infty]$. For $p>2$
it is easy to check that $\sum_{i=1}^k\|y_i\| \leq k^{1/p'} \|\sum_{i=1}^k
y_i\|$ for any disjointly
supported
$y_1,\ldots,y_k$ on the canonical basis of $Y_p$ (see e.g. \cite[Lemma~2.4]{FG1}). Therefore by \cite[Lemma~9.3]{K1}
$Y_p$ satisfies a lower
$r$-estimate for any $r>p$ and therefore has cotype $r$ \cite{LT2} p.88. On the
other hand $Y_p$ contains
$\ell_p^n$'s  uniformly and therefore does not have cotype $q$ for $q<p$.
Therefore the spaces $Y_p$, $p>2$,
are mutually nonisomorphic.
\end{proof}

\end{document}